\documentclass{elsart}
\usepackage{amsmath}
\usepackage{amsfonts}
\usepackage{enumerate}
\usepackage[numbers,sort&compress]{natbib}

\newcommand{\reals}{\mathbb R}          

\theoremstyle{plain}
\newtheorem{theorem}{Theorem}
\newtheorem{lemma}{Lemma}
\newtheorem{proposition}[lemma]{Proposition}

\theoremstyle{definition}

\newenvironment{proof}[1][]{\textit{Proof#1}. }{\qed\par}
\newtheorem{remark}{Remark}

\newcommand{\Laplacian}{L}
\newcommand{\ac}{\alpha}

\begin{document}


\begin{frontmatter}

\title{Algebraic Connectivity and Degree Sequences of Trees}

\author[IU]{T\"urker B{\i}y{\i}ko\u{g}lu} and
\ead{turker.biyikoglu@isikun.edu.tr}
\ead[url]{http://math.isikun.edu.tr/turker/}
\author[WU]{Josef Leydold\corauthref{cor}}
\ead{Josef.Leydold@statistik.wu-wien.ac.at}
\ead[url]{http://statistik.wu-wien.ac.at/\~{}leydold/}

\address[IU]{Department of Mathematics,
  I\c{s}{\i}k University,
  \c{S}ile 34980, Istanbul, Turkey}
\address[WU]{Department of Statistics and Mathematics,
  University of Economics and Business Administration,
  Augasse 2-6, A-1090 Wien, Austria}

\corauth[cor]{Corresponding author. Tel +43 1 313 36--4695. FAX +43 1 313 36--738}



\begin{abstract}
  We investigate the structure of trees that have minimal algebraic
  connectivity among all trees with a given degree sequence.
  We show that such trees are caterpillars and that the vertex degrees
  are non-decreasing on every path on non-pendant vertices starting at
  the characteristic set of the Fiedler vector.
\end{abstract}

\end{frontmatter}


\markboth{T.~B{\i}y{\i}ko\u{g}lu and J.~Leydold}{%
  Algebraic Connectivity and Degree Sequences of Trees}


\section{Introduction}

Let $G(V,E)$ be a simple (finite) undirected graph with vertex set $V$
and edge set $E$. 
The total number of vertices is denoted by $n$.
The \emph{Laplacian} of $G$ is the matrix 
\begin{equation}
  \Laplacian(G) = D(G) - A(G)\;,  
\end{equation}
where $A(G)$ denotes the adjacency matrix of the graph and $D(G)$ is
the diagonal matrix whose entries are the vertex degrees, i.e.,
$D_{vv} = d(v)$, where $d(v)$ denotes the degree of vertex $v$.
We write $\Laplacian$ for short if there is no risk of confusion.
For graphs with weights $w(e)$ for each edge $e\in E$ the Laplacian
is defined analogously where the adjacency matrix contains the edge
weights and the diagonal entries of $D(G)$ are the sums of the
weights of the edges at the vertices of $G$, i.e.\
$D_{vv} = \sum_{e=vu\in E} w(e)$.

The Laplacian $\Laplacian$ is symmetric and all its eigenvalues
are non-negative. The first eigenvalue is always $0$.
The second smallest eigenvalue, denoted by $\ac(G)$ in the following,
has become quite popular and is called the \emph{algebraic
  connectivity} of $G$ by \citet{Fiedler:1973a}.  
It allows some conclusions about the connectedness of the graph. A
graph $G$ is connected if and only if $\ac(G)\not=0$. Moreover,
$\ac(G)$ is a lower bound for the vertex and edge connectivities of $G$.  
Hence properties of the algebraic connectivity has been investigated
in the literature. In particular many upper and lower bounds have been
shown. We refer to the recent survey by \citet{Abreu:2007a} and the
references cited therein.
Other authors have ordered trees by their algebraic connectivity
\citep{Shao;Guo;Shan:2008a} or 
characterize extremal graphs, i.e., graphs that have
minimal algebraic connectivity among all graph within particular graph
class. \citet{Godsil;Royle:2001a} assume that \emph{graphs with small values
of $\ac(G)$ tend to be elongated graphs of large diameter with bridges.}
For example, for trees on $n$ vertices with a fixed diameter
the algebraic connectivity is minimized for paths with stars of
(almost) equal size attached to both ends, see
\citep{Fallat;Kirkland:1998a}. Cubic graphs with minimal algebraic
connectivity look like a ``string of pearls'', see
\citep{Brand;Guiduli;Imrich:2007a}. 
\citet{Belhaiza;etal:2005a} used the AGX-system which raised the
conjecture that the connected graphs $G\not=K_n$ with minimal
algebraic connectivity are all so called $(n,p,t)$-path-complete
graphs.

In this note we are interested in the structure of trees which have
minimal algebraic connectivity among all trees with a given degree
sequence. 
(Recall that a sequence $\pi=(d_0,\ldots,d_{n-1})$ of non-negative
integers is called \emph{degree sequence} if there exists a graph $G$
for which $d_0,\ldots,d_{n-1}$ are the degrees of its vertices. 
We refer the reader to \citep{Melnikov;etal:1994a} for
relevant background on degree sequences.)
We call a degree sequence for a tree a \emph{tree sequence}.
We show that such a tree is a \emph{caterpillar}, i.e., a tree
in which the removal of all pendant vertices (vertices of degree 1)
gives a path. 
For further characterization we use a result of \citet{Fiedler:1975b}
about eigenvectors of the second smallest eigenvalue which are called
\emph{Fiedler vectors}:
The subgraph induced by non-positive vertices of any Fiedler vector
(i.e., vertices with non-positive valuation) and the subgraph induced
non-negative vertices are both connected. 
Such connected subgraphs are called \emph{weak nodal domains}
\citep{Davies;etal:2001a,Biyikoglu;Leydold;Stadler:2007a} (in analogy
to eigenfunctions of the Laplace-Beltrami operator on manifolds
\citep{Courant;Hilbert:1953a,Cheng:1976a}),
or \emph{Perron components} \citep{Kirkland;Neumann;Shader:1996a}. 
The two nodal domains of a Fiedler vector of a tree 
are either separated by an edge (\emph{characteristic edge}) or there
is a single vertex (\emph{characteristic vertex}) where the Fiedler
vector vanishes \citep{Fiedler:1975b}. 
On each of these nodal domains we can declare a Dirichlet matrix whose
first eigenvalue is exactly the algebraic connectivity of the original
graph. Thus we arrive at the following necessary condition.

\clearpage

\begin{theorem}
  \label{thm:1}
  Let $T$ be a tree that has minimal algebraic connectivity among all
  trees with given degree sequence $\pi=(d_0,\ldots,d_{n-1})$.
  Then $T$ is a caterpillar. 
  Moreover, if $P$ is the path induced by all non-pendant vertices of
  $T$ with non-negative (non-positive) valuation, then its degree
  sequence is monotone with a minimum at the characteristic vertex or
  edge.
\end{theorem}

\begin{rem}
  It is an open problem how the degree sequence $\pi$ has to be
  partitioned for the two nodal domains to obtain a tree with minimum 
  algebraic connectivity. We ran some computational experiments but
  could not detect any general pattern.
\end{rem}

For the proof of this theorem we use the concept of \emph{geometric
  nodal domains} and \emph{Dirichlet matrix} introduced in 
\citep{Biyikoglu;Leydold;Stadler:2007a}. This concept is described in 
Sect.~\ref{sec:nodal-domain}. We then can use perturbation of trees to
obtain results for the first Dirichlet eigenvalue for each of the two
nodal domains of the Fiedler vector
(Sect.~\ref{sec:rooted-tree}) which are then used to proof the
theorem in Sect~\ref{sec:proof}.
Our approach is related to the concept of \emph{Perron components} and 
\emph{``bottleneck'' matrix} introduced in
\citep{Kirkland;Neumann;Shader:1996a}.
Thus it could also be used to verify the results of 
\citep{Kirkland;Neumann:1997a}
 (e.g., their Thm.~5 can be deduced from
 Lemma~\ref{lem:dirichlet-perturbation} below).


\section{Nodal Domains and Dirichlet Matrix}
\label{sec:nodal-domain}

A \emph{graph with boundary} $G(V_0\cup\partial V,E_0\cup\partial E)$ 
consists of a (non-empty) set of interior vertices $V_0$, boundary
vertices $\partial V$, interior edges $E_0$ that connect interior
vertices, and boundary edges $\partial E$ that join interior vertices
with boundary vertices. There are no edges between two boundary
vertices. 
The \emph{Dirichlet matrix} $\Laplacian_0$ is
the matrix obtained from the graph Laplacian $\Laplacian$ by deleting
all rows and columns that correspond to boundary vertices. 
This definition is motivated by the concept of \emph{geometric
  realization} of a graph, see
\citep{Friedman:1993a,Biyikoglu;Leydold;Stadler:2007a}.
The first Dirichlet eigenvalue $\nu(G)$ is strictly positive.
If the graph induced by the interior vertices is connected, then
$\nu(G)$ is simple and there exists an eigenvector which is strictly
positive at all interior vertices.

When a Fiedler vector $f$ of a tree has a characteristic vertex $v_0$
then each of the two weak nodal domains can be seen as a graph with
$v_0$ as its boundary vertex. Then the algebraic connectivity $\ac(G)$
is exactly the first Dirichlet eigenvalue of each nodal domain, with
the Fiedler vector restricted to the respective interior vertices as
their eigenvectors (see also \citep{Bapat;Pati:1998a}).
In the other case when the two nodal domains of the Fiedler
vector are separated by a characteristic edge $e=uw$, we split
this edge of weight $1$ into edges $e_1=uv_0$ and $e_2=v_0w$ with
weights $w_1=|f(w)-f(u)|/|f(u)|$ and $w_2=|f(w)-f(u)|/|f(w)|$
by inserting a new vertex $v_0$. 
(In the geometric realization of a graph edges of weight $w$
correspond to arcs of length $1/w$.)
By this procedure the algebraic connectivity remains unchanged and
$v_0$ becomes the characteristic vertex of the Fiedler vector of the
new graph \citep[Lemma~3.14]{Biyikoglu;Leydold;Stadler:2007a}.
In either case we construct two graphs with boundary whose first
Dirichlet eigenvalues coincide with $\ac(G)$. We call these graphs the
\emph{geometric nodal domains} of $G$.
Thus we can prove our theorem by looking at the first Dirichlet
eigenvalue of its nodal domains.

\begin{remark}
  The concept of geometric nodal domains is defined analogously for
  arbitrary eigenfunctions of connected graphs.
\end{remark}

The Rayleigh quotient associated to the Laplacian matrix $L$ is
defined by
\begin{equation}
  \mathcal{R}_\Laplacian(f)
  = \frac{\langle f, \Laplacian f\rangle}{\langle f,f \rangle}
  = \frac{\sum_{uv\in E} w(uv)(f(u)-f(v))^2}{\sum_{v\in V} f(v)^2}\;.
\end{equation}

The following result characterizes the first Dirichlet eigenvalue
$\nu(G)$ and the algebraic connectivity $\ac(G)$ of some graph $G$. It
immediately follows from the Courant-Fisher Theorem.
\begin{proposition}
  \label{prop:minimax}
  For a graph with boundary $G(V_0\cup\partial V, E_0\cup\partial E)$ 
  we have
  \begin{equation}
    \nu(G) = \min_{f\in\reals^{|V_0|},\,f\not=0} 
    \mathcal{R}_{\Laplacian_0}(f)
  \end{equation}
  Moreover, $f\not=0$ is an eigenvector of the first Dirichlet
  eigenvalue $\nu(G)$ of $\Laplacian_0$ if and only if
  $\mathcal{R}_{\Laplacian_0}(f)=\nu(G)$.

  For a graph $G(V,E)$ we have
  \begin{equation}
    \ac(G) = \min_{f\in\reals^{|V|},\,f\not=0,\,\sum f(v)=0} 
    \mathcal{R}_{\Laplacian}(f)
  \end{equation}
  Moreover, $f\not=0$ is an eigenvector of the second
  Laplacian eigenvalue $\ac(G)$ (i.e.\ a Fiedler vector) if and only
  if $\sum_{v\in V} f(v)=0$ and $\mathcal{R}_\Laplacian(f)=\ac(G)$.
\end{proposition}


\section{First Dirichlet Eigenvalues of Rooted Trees}
\label{sec:rooted-tree}

Geometric nodal domains of Fiedler vectors of trees are rooted trees
$T_0$ where its root $v_0$ is its only boundary vertex. One of its 
boundary edges has weight $w_0$ with $1/w_0\in(0,1]$ whereas all other
(boundary and interior) edges have weight $1$.
The following lemma immediately follows from
\citep[Thm.~(3,14)]{Fiedler:1975b}. 

\begin{lemma}
  \label{lem:increasing}
  Let $T_0$ be a tree with a single boundary vertex $v_0$ and 
  $f$ a non-negative eigenvector corresponding to the first
  Dirichlet eigenvalue $\nu(T_0)$. Then on very every simple path
  starting at $v_0$, $f$ is either strictly increasing or constant
  zero.
\end{lemma}

A \emph{branch at vertex $u$} of a tree with \emph{root} $v_0$ is a
maximal subtree of $G\setminus\{u\}$ that does not contain $v_0$.

\begin{lemma}
  \label{lem:dirichlet_necessary}
  If a tree $T_0$ has minimal first Dirichlet eigenvalue among all
  rooted trees with given degree sequence,
  then $T_0$ is a caterpillar where at most one neighbors of its
  root is not a pendant vertex.
\end{lemma}
\begin{proof}
  Let $T_0$ have minimal first Dirichlet eigenvalue, i.e.,
  $\nu(T_0)\leq\nu(T'_0)$ for all rooted trees $T'_0$ with the same
  degree sequence. 
  Assume first that there is only one boundary edge $xv_0$, then
  $\nu(T_0)$ is simple and there exists an eigenvector $f$ with
  $f(v_0)=0$ and $f(u)>0$ for all $u\not=v_0$.
  Now suppose $T_0$ is not a caterpillar with the proposed
  property. Then there exist two simple paths
  $(v_0,\ldots,v_{i-1},x_i,\ldots,x_j)$ and
  $(v_0,\ldots,v_{i-1},y_i,\ldots,y_k)$ where $x_i\not= y_i$,
  $j,k>i$, and where $x_j$ and $y_k$ are pendant vertices.
  By Lemma~\ref{lem:increasing}, $f$ is strictly increasing on each of
  these. Without loss of generality we assume
  $f(x_j)>f(y_i)$. Otherwise we have 
  $f(x_i)<f(x_j)\leq f(y_i)<f(y_k)$ and we exchange the role of the
  two paths.
  Now we construct a new graph $T'_0$ by rearranging edges in $T_0$.
  Replace all edges $y_it$ where $t\not=v_{i-1}$ by edges $x_jt$.
  Notice that this rearrangement does not change the degree sequence.
  We construct a new vector $f'$ such that
  $f'(u)=\max(f(u),f(x_j))$ for all vertices $u$ that are in a branch
  at $y_i$ in $T_0$, and $f'(u)=f(u)$ for all others. Notice that
  $\sum_{v\in V} f'(v)^2 \geq \sum_{v\in V} f(v)^2$ and
  $\sum_{uv\in E'} w(uv)(f'(u)-f'(v))^2 < 
  \sum_{uv\in E} w(uv)(f(u)-f(v))^2$ and thus
  $\nu(T'_0)\leq\mathcal{R}_\Laplacian(f')
  <\mathcal{R}_\Laplacian(f)=\nu(T_0)$, a contradiction to our
  assumption that $T_0$ has minimal first Dirichlet eigenvalue.
  \\
  For the case where $T_0$ has two or more boundary edges then
  $T_0\setminus\{v_0\}$ consists of several branches at root $v_0$.
  There exists an eigenvector $f$ that is positive on exactly one of
  these and zero on all others. 
  Then all these other branches must be pendant vertices, since
  otherwise we could apply the same rearrangement of edges and obtain
  a tree $T'_0$ with the same degree sequence and strictly smaller first
  Dirichlet eigenvalue.
\end{proof}

The \emph{trunk} of a rooted caterpillar $T_0$ is a longest path
starting at root $v_0$. Notice that the trunk is terminated by $v_0$
and a pendant vertex (its \emph{head}).
Let $h(v)$ denote the geodetic distance between vertex $v$ and root
$v_0$ (\emph{height}).
Now construct a new rooted graph $T_0'$ by one of the
following perturbations:
\begin{enumerate}[(P1)]
\item
  Replace edge $wv_i$ by an edge $wv_j$, where $w$ ($\not=v_0$) is a
  pendant vertex, and $v_i$ and $v_j$ are trunk vertices with
  $h(v_i)<h(v_j)$;
\item 
  insert a vertex $w$ and add a new edge $wv_j$ to one of the trunk
  vertices $v_j$ ($\not=v_0$).
\end{enumerate}
Notice that in both cases the trunk of $T_0'$ is longer than that of
$T_0$ if $v_j$ is the head of the trunk of $T_0$.

\begin{lemma}
  \label{lem:dirichlet-perturbation}
  Let $T_0$ a rooted tree and construct a new tree $T_0'$ as described
  above. Then $\nu(T'_0) < \nu(T_0)$.
\end{lemma}
\begin{proof}
  Let $f$ be a nonnegative eigenfunction to the first Dirichlet
  eigenvalue $\nu(T_0)$. 
  For case (P1) we construct function $f'$ on $T_0'$ by
  $f'(v)=f(v)$ for all $v\not=w$ and $f'(w)=\max(f(w),f(v_j))$.
  By Lemma~\ref{lem:increasing}, $f$ is strictly increasing on the
  trunk of $T$. Thus we find analogously to the proof of
  Lemma~\ref{lem:dirichlet_necessary},
  $\nu(T'_0)\leq\mathcal{R}_\Laplacian(f')
  <\mathcal{R}_\Laplacian(f)=\nu(T_0)$ as proposed.
  For case (P2) we set $f'(v)=f(v)$ for all $v\not=w$ and $f'(w)=f(v_j)>0$
  and thus $\sum_{v\in V} f'(v)^2 > \sum_{v\in V} f(v)^2$ and the
  statement follows.
\end{proof}

\begin{lemma}
  \label{lem:dirichlet_necessary_sufficient}
  A tree $T_0$ has minimal first Dirichlet eigenvalue among all rooted
  trees with given degree sequence if and only if $T_0$ is a caterpillar
  where the degrees are non-decreasing on the path starting at
  root $v_0$ and induced by all non-pendant vertices of $T_0$.
\end{lemma}
\begin{proof}
  By Lemma~\ref{lem:dirichlet_necessary}, $T_0$ is a caterpillar.
  Let $(v_0,v_1,\ldots,v_k,v_{k+1})$ be the trunk of $T_0$.
  If $T_0$ does not have increasing degrees, then there exist vertices
  $v_i$ and $v_j$, $i<j$, in this path with $2\leq d(v_j)<d(v_i)$.
  Thus we can replace $c=d(v_j)-d(v_i)$ edges and get a new graph
  $T'_0$ with the same degree sequence as $T_0$. By 
  Lemma~\ref{lem:dirichlet-perturbation}, $\nu(T_0) < \nu(T'_0)$, a
  contradiction to our assumptions.
\end{proof}


\section{Proof of the Theorems}
\label{sec:proof}

We now show our result by gluing two rooted trees together.

\begin{lemma}
  \label{lem:union}
  Let $T_1$ and $T_2$ be two trees with one boundary vertex.
  Construct a new tree $T$ without boundary by identifying the 
  boundary vertices of these trees and turning the new vertex into an
  interior vertex. \\
  Then $\ac(T)\leq\max(\nu(T_1),\nu(T_2))$.
  The inequality is strict if $\nu(T_1)\not=\nu(T_2)$.
\end{lemma}
\begin{proof}
  We assume without loss of generality that $\nu(T_1)\geq\nu(T_2)$.
  Let $V_1$ and $V_2$ denote the respective vertex sets and 
  $f_1$ and $f_2$ be corresponding non-negative Dirichlet
  eigenvectors such that 
  $\sum_{v\in V_1} f_1(v) = \sum_{w\in V_2} f_2(v)$.
  Construct a vector $f$ on $T$ by $f(v)=f_1(v)$ for all $v\in V_1$
  and $f(u)=-f_2(u)$ for $u\in V_2$.
  Notice that for positive numbers $x,y,a,b>0$ we find
  $\frac{x+a}{y+b} \leq \frac{x}{y}$ if and only if 
  $\frac{a}{b} \leq \frac{x}{y}$ and that either both or none of the
  equalities hold.
  Then we find by Proposition~\ref{prop:minimax}
  \begin{eqnarray*}
    \ac(T) &=& \min_{\sum g(v)=0} \mathcal{R}_{\Laplacian}(g) \\
    &\leq& \mathcal{R}_{\Laplacian}(f) 
    = \frac{\sum_{uv\in E} w(uv)(f(u)-f(v))^2}{\sum_{v\in V} f(v)^2} \\
    &=& \frac{\sum_{uv\in E_1} w(uv)(f_1(u)-f_1(v))^2
      + \sum_{uv\in E_2} w(uv)(f_2(u)-f_2(v))^2}{
      \sum_{v\in V_1} f_1(v)^2 + \sum_{v\in V_2} f_2(v)^2} \\
    &\leq& \frac{\sum_{uv\in E_1}
      w(uv)(f_1(u)-f_1(v))^2}{\sum_{v\in V_1} f_1(v)^2}
    = \nu(T_1)\;.
  \end{eqnarray*}
  Moreover, $\ac(T)<\nu(T_1)$ whenever $\nu(T_1)>\nu(T_2)$ and thus
  the second statement follows.
\end{proof}

\begin{proof}[ of Theorem~\ref{thm:1}]
  Assume $T$ has minimal algebraic connectivity $\ac(T)$ and let
  $T_1$ and $T_2$ be its the two geometrical nodal domains.
  Then $\nu(T_1)=\nu(T_2)=\ac(T)$. Both subtrees must be
  caterpillars as described in
  Lemma~\ref{lem:dirichlet_necessary_sufficient}. Otherwise, if (say)
  $T_1$ does not have this property we could replace it by a tree
  $T_1'$ with the same corresponding edge weights as in $T_1$ and with
  with $\nu(T_1')<\ac(T)$. Consequently we could construct a new tree
  graph $T'$ with the same degree sequence as in $T$ but with
  $\ac(T')<\ac(T)$ by Lemma~\ref{lem:union}, a contradiction.
  Thus the statement follows.
\end{proof}



\end{document}